\documentclass[12pt]{amsart}
\newtheorem{thm}{Theorem}

\newtheorem*{problem*}{Problem}
\newcommand{\Z}{\mathbb Z}
\newcommand{\N}{\mathbb N}
\newcommand{\R}{\mathbb R}

\begin{document}
\title[Relative Variational Principle: A survey]{The Relative Variational Principle by Ledrappier and Walters: A survey}
\author{Anthony Quas}
\address{Department of Mathematics and Statistics, University of Victoria, Victoria BC, CANADA V8W 3R4}
\subjclass[2020]{37D35}
\maketitle
\begin{abstract}Ledrappier and Walters's article ``A Relativised Variational Principle for Continuous Transformations", 
J. Lond. Math. Soc. (2) 16 (1977), no.3, 568-576) is a landmark 
in the development of Thermodynamic Formalism. This survey, aimed at newcomers to the field and experts in adjacent 
fields discusses the background, the Ledrappier-Walters article and some subsequent developments in the field.
\end{abstract}

\section{Introduction}
Entropy and pressure are a key part of a remarkable family of ideas, known as \emph{Thermodynamic Formalism}.
The ideas, now ubiquitous in ergodic theory and dynamical systems, were originally based on 
introducing ideas from 19th Century thermodynamics into modern dynamical systems. 
The early stages of the development of the thermodynamic formalism
took place over two decades, from the late 1950's to the late 1970's, with
 Kolmogorov and Sinai introducing the notion of measure-theoretic entropy in the late 1950's; 
 Adler, Konheim and McAndrew introducing topological entropy in the mid 1960's, a new definition of topological entropy, due to Bowen,
in the early 1970's, the Variational Principle for topological entropy (due to Goodman, with contributions from
Goodwyn and Dinaburg) in the early 1970's and the Variational Principle for pressure by Ruelle
and Walters in the mid 1970's. 

Thermodynamic formalism now plays a key role
in applications of dynamical systems to number theory, information theory, geometry and other fields. 
The paper, ``A relativised Variational Principle for continuous transformations" by Ledrappier and Walters is one of the key contributions to the 
theory. Indeed, the theory is now so well known and easily applied that someone
entering the field may not recognize how revolutionary
this circle of ideas was. 

In this survey article, we aim to present an overview of the field for newcomers, as well as researchers
in adjacent fields. We mostly work in the setting of the article \cite{LedrappierWalters} of Ledrappier and Walters, 
and restrict our attention to continuous dynamical systems defined by a map $T\colon X\to X$, 
where $X$ is a compact metric space, even though
many of the results have versions that hold valid in a wider context. In many places, we restrict even further, to the case where
$X$ is a \emph{subshift}, a closed shift-invariant subset of a set $A^\Z$ (equipped with the product topology),
 where $A$ is a finite set, and $T$ is the shift map, 
$T(x)_n=x_{n+1}$. While the specific metric is unimportant in the case of compact $X$, 
on subshifts we always use the metric
$d(x,y)=2^{-n(x,y)}$ where $n(x,y)=\min\{|n|:x_n\ne y_n\}$.

I thank the referees for helpful suggestions for improving this article.

\section{Background}
\subsection{Measure-theoretic entropy}
If $(X,\mu)$ is a probability space, the entropy of a countable partition, 
$\mathcal P=\{A_1,A_2,\ldots\}$ is defined by $H_\mu(\mathcal P)=-\sum_i \mu(A_i)\log\mu(A_i)$,
where we adopt the convention $0\log 0=0$. 
For two countable partitions, $\mathcal P$ and $\mathcal Q$, their \emph{join}, $\mathcal P\vee\mathcal Q$, is 
$\{A\cap B\colon A\in\mathcal P, B\in\mathcal Q\}$. 
If $T$ is a measure-preserving transformation of $(X,\mu)$, the \emph{measure-theoretic entropy of $T$ with respect to $\mathcal P$} is defined by
$$
h_\mu(T,\mathcal P)=\lim_{n\to\infty}\frac 1nH\Big(\bigvee_{i=0}^{n-1}T^{-i}\mathcal P\Big).
$$
The existence of the limit follows from the measure-preserving property, the fact that 
$H_\mu(\mathcal P\vee\mathcal Q)\le H_\mu(\mathcal P)+H_\mu(\mathcal Q)$
and Fekete's lemma.
Finally, the \emph{measure-theoretic entropy of $T$} is defined by
$$
h_\mu(T)=\sup_{\mathcal P} h_\mu(T,\mathcal P),
$$
where the supremum can be taken either over finite partitions, or partitions $\mathcal P$ satisfying $H_\mu(\mathcal P)<\infty$. While existence is
clear, it is remarkable that $h_\mu(T)$ is often positive and finite. 

Surprisingly, it is often relatively easy to calculate: $\mathcal P$ is called a 
\emph{generating
partition} if there is a set $X_0$ of measure 1 such that if $x,y$ are distinct elements of $X_0$, then there exists $n$ such that
$T^n(x)$ and $T^n(y)$ lie in distinct elements of $\mathcal P$.
(If $T$ is invertible, $n$ runs over $\Z$; while if $T$ is not invertible, $n$ runs over $\N_0$).
 In the case that $\mathcal P$ is a generating partition,
the supremum in the definition of the measure-theoretic entropy is attained: $h_\mu(T)=h_\mu(T,\mathcal P)$. 

If $X$ is a subshift, a simple example
of a generating partition is $\mathcal P=\{[a]\colon a\in A\}$, where
$[a]$ is the \emph{1-cylinder set} $\{x\in X\colon x_0=a\}$. 

Since measure-theoretic entropy is invariant under isomorphism 
(essentially measurable changes of coordinates), it provides a mechanism
for showing that measure-preserving transformations fail to be isomorphic. Indeed,
Kolmogorov and Sinai \cite{KolmogorovSinai} introduced measure-theoretic entropy to show that the Bernoulli shifts with 2 symbols with equal weight; 
and with 3 symbols with equal weight
are not measure-theoretically isomorphic (the first has entropy $\log 2$ and the other $\log 3$). 

\subsection{Shannon-MacMillan-Breiman}
The Shannon-MacMillan-Breiman theorem \cite{Breiman} gives a very natural interpretation of the measure-theoretic entropy. In the case
where  $T\colon X\to X$
is a subshift and $\mu$ is an ergodic invariant measure, then
$$
-\tfrac 1n\log\mu([x_0\ldots x_{n-1}])\to h_\mu(X)\text{\quad for $\mu$-a.e.\ $x$},
$$
where $[x_0\ldots x_{n-1}]$ is the \emph{$n$-cylinder set}
$\{y\colon y_i=x_i\text{ for $0\le i<n$}\}$.

\subsection{Topological entropy}\label{sec:topent}
For a continuous dynamical system, $T\colon X\to X$, we now introduce the \emph{topological entropy}.
This was defined
by Adler, Konheim and McAndrew \cite{AKMcA}, who also conjectured the relationship with 
measure-theoretic entropy that we describe
in Section \ref{sec:vpent}.
Alternative formulations were given by Dinaburg \cite{Dinaburg} and Bowen \cite{Bowen:EntropyforGroupEndomorphisms}. 
We now describe Bowen's formulation.
We view topological entropy as a measure of the exponential growth rate of the number of orbit segments 
of length $n$ as a function of the length. Since there are, in all non-trivial cases,
an infinite number of orbits, one counts the number of ``distinguishable" orbit segments of length $n$, 
where two orbit segments are
distinguishable if they diverge by at least a pre-determined threshold by time $n$ (i.e. the $n$-orbits of 
$x$ and $y$ are $\epsilon$-distinguishable if $d(T^ix,T^iy)\ge \epsilon$ for some $0\le i<n$). 

Formally, the topological entropy is given by
$$
h_\text{top}(T)=\lim_{\epsilon\to 0}
\limsup_{n\to\infty}\tfrac 1n\log N_\epsilon(n),
$$
where $N_\epsilon(n)$ is the maximal cardinality of an $(n,\epsilon)$-separated subset $S$
of $X$: sets such that for any distinct $x,y\in S$, $d(T^ix,T^iy)\ge \epsilon$
for some $0\le i<n$. Compactness of $X$ together with continuity of $T$ implies 
$N_n(\epsilon)$ is finite for all $n$ and $\epsilon$.
The combinatorial bound
$N_{m+k}(\epsilon)\le N_k(\frac \epsilon 2)N_m(\epsilon)$ can be used to show 
$\limsup_{n\to\infty}\frac 1n\log N_n(\epsilon)
\le \frac 1k\log N_k(\frac\epsilon 2)$ for any $k$, so that the inner limit is finite for any $\epsilon>0$. 
We observe that $\limsup_{n\to\infty}\frac 1n\log N_n(\epsilon)$ is non-decreasing as $\epsilon$ is shrunk to zero. 
The limit may, however, diverge. 

A refinement of the above definition of topological entropy that we require later is a measure
of the complexity \emph{in a set $S$}:
$$
h_\text{top}(T,S)=\lim_{\epsilon\to 0}
\limsup_{n\to\infty}\tfrac 1n\log N_\epsilon(n,S),
$$
where $N_\epsilon(n,S)$ is the maximal cardinality of an $(n,\epsilon)$-separated subset of $S$.
We emphasize that $S$ is not required to be an invariant set. 

For subshifts with a finite alphabet, the topological entropy is always finite, and is often easy to compute.
In this case, one can check that a maximal $(n,1)$-separated set consists of an arbitrary set of points $S$,
one from each (non-empty) $n$-cylinder set.
Similarly, maximal $(n,2^{-k})$-separated sets are also arbitrary representatives of cylinder sets of the form 
$\{x\colon x_i=a_i\text{ for $i=-k.\ldots,n+k-1$}\}$. Hence for subshifts, maximal 
$(n,2^{-k})$-separated sets are the images under 
$T^k$ of maximal $(n+2k,1)$-separated sets. 
From the above, we obtain in the case of a subshift
$$
h_\text{top}(T)=\lim_{n\to\infty}\tfrac 1n\log C_n(X),
$$
where $C_n(X)$ is the number of non-empty $n$-cylinder sets.

\subsection{The Variational Principle for topological entropy}\label{sec:vpent}
As conjectured by Adler, Konheim and McAndrew, if $(X,d)$ is compact and
$T:X\to X$ is a continuous dynamical system, the topological entropy is 
related to the the measure-theoretic entropy of $T$ by the formula
$$
h_\text{top}(T)=\sup_\mu h_\mu(T),
$$
where the supremum is taken over the collection of $T$-invariant Borel probability measures. This was proved, 
in varying degrees of generality,
by Dinaburg \cite{Dinaburg}, Goodwyn \cite{Goodwyn}, Goodman \cite{Goodman}. In the case where $T\colon X\to X$
 is a subshift (and, more generally,
if $T$ satisfies the \emph{expansiveness} condition), the supremum is attained: 
there is at least one ergodic \emph{measure of maximal entropy}. That is, a measure $\mu$ on $X$ 
such that $h_\mu(T)=h_\text{top}(T)$. The underlying reason is that if $T$ is expansive,
the map $\mu\mapsto h_\mu(T)$ is upper semi-continuous with respect to the weak$^*$-topology on the
compact set  of $T$-invariant measures.
The existence and/or multiplicity of ergodic measures of maximal entropy 
in general dynamical systems
has received a lot of attention. A system is said to be 
\emph{intrinsically ergodic} if it has a unique measure of maximal entropy. A well-known
class of systems with this property
is the class of irreducible shifts of finite type \cite{ParryIntrinsic}.

If $\mu$ is an ergodic measure of maximal entropy, a consequence of the
Shannon-MacMillan-Breiman theorem is
 for $\mu$-a.e.\ $x$, for all $\epsilon>0$,
\begin{equation}\label{eq:SMB}
e^{-n\epsilon}\le \frac{\mu([x_0\ldots x_{n-1}])}{1/C_n(X)}\le e^{n\epsilon}
\end{equation}
for all sufficiently large $n$. That is, for an ergodic measure of maximal entropy,
measures of most cylinder sets agree, up to a sub-exponential factor, with $1/C_n(X)$.

As a critical caveat to the above, the word \emph{most} here refers to most cylinder sets
\emph{with respect to $\mu$} and not with respect to counting. That is, for any $\epsilon>0$ and
$\delta>0$, for all sufficiently large $n$, there is a collection of $n$-cylinder sets of 
combined measure at least $1-\delta$, each satisfying \eqref{eq:SMB}.

\subsection{Pressure}

An important generalization of topological entropy, called \emph{(topological) pressure}, was 
introduced by Ruelle \cite{RuelleRigorous,RuelleVP}. Ruelle's formulation used ideas from statistical
physics , as well as Bowen's definition of topological entropy. Ruelle's topological 
pressure was for dynamical systems satisfying 
the additional hypotheses of expansiveness and specification, which we will not describe here, but which 
imply that pressure and topological entropy can be computed using the sets of period $n$ points for
each $n$ instead of
$(n,\epsilon)$-separated sets. Additionally, in keeping with the statistical 
mechanical background, Ruelle's definition allowed for actions of $\Z^d$ rather than just $\Z$. That is, Ruelle considered
the case of the commuting action of $d$ continuous transformations. The prototype for this 
is the $d$ coordinate translation maps
of $\Z^d$ configurations of a system of particles. 

Ruelle's definition was subsequently adapted by Walters \cite{WaltersVP} to cover general continuous
dynamical systems (actions of $\Z$ rather than $\Z^d$) without assuming either 
expansiveness or specification. We now review Walters's definition,
which is closely related to Bowen's definition of topological entropy, described above, in terms of $(n,\epsilon)$-separated
sets. 
Given a continuous function $\phi$, we define $P_n(T,\phi,\epsilon)$ by
$$
P_n(T,\phi,\epsilon)=\sup_{S\text{ $(n,\epsilon)$-separated}}\sum_{x\in S}\exp\left(\sum_{j=0}^{n-1}\phi(T^nx)\right).
$$
The pressure of $\phi$ with respect to $T$, $P_T(\phi)$ is then defined by
$$
P_T(\phi)=\lim_{\epsilon\to 0}\limsup_{n\to\infty}\tfrac 1n\log P_n(T,\phi,\epsilon).
$$
Since it is straightforward to see that $P_n(T,\phi,\epsilon)\ge P_n(T,\phi,\epsilon')$ if $\epsilon<\epsilon'$, 
existence of the limit is clear. 
Also, if $\phi$ is the 0 function, the pressure agrees with Bowen's definition of the topological entropy above. 
One can verify that $|P_T(\phi)-h_\text{top}(T)|\le \|\phi\|$, so that for continuous functions on compact spaces, the pressure
is finite if and only if the topological entropy is finite. 

As with topological entropy, one may define the pressure on a subset by
$$
P_T(\phi,A)=\lim_{\epsilon\to 0}\limsup_{n\to\infty}\tfrac 1n\log P_n(T,\phi,\epsilon,A),
$$
where $P_n(T,\phi,\epsilon,A)$ is the supremum over $(n,\epsilon)$-separated subsets of $A$.

In the case of subshifts, similarly to the situation for topological entropy, $\limsup \frac 1n P_n(T,\phi,\epsilon)$
is the same for all $\epsilon\le 1$, so we obtain the simpler formula

$$
P_T(\phi)=\limsup_{n\to\infty} \tfrac 1n\log P_n(T,\phi,1).
$$

For subshifts, using uniform continuity of $\phi$, one may show that
$$
\frac 1n\log\left[\left(\sup_{A\in \mathcal R_n(X)} \sum_{x\in A}\exp S_n\phi(x)\right)/
\left(\inf_{A\in \mathcal R_n(X)} \sum_{x\in A}\exp S_n\phi(x)\right)\right]\to 0,
$$
where $\mathcal R_n(X)$ denotes the collection of all sets of size $C_n(X)$ with one representative from each 
non-empty cylinder set, and $S_n\phi(x)=\sum_{i=0}^{n-1}\phi(T^ix)$.
Hence, for subshifts
\begin{equation}\label{eq:Pcalc}
P_T(\phi)=\limsup_{n\to\infty}\frac 1n\log \sum_{x\in A_n}\exp S_n\phi(x),
\end{equation}
where $A_n$ is an arbitrary element of $\mathcal R_n(X)$ for each $n$. In a final simplification, the sequence
actually converges, so the limit superior can be replaced by a limit.

\subsection{The Variational Principle for pressure}

Walters's paper above (and Ruelle's paper in the context of expansive $\Z^d$ actions with specification)
proved a Variational Principle for Pressure. We state Walters's version here:

\begin{thm}[Variational Principle]
Let $X$ be a compact metric space and let $T\colon X\to X$ be a continuous dynamical system. Let $\phi\colon X\to\R$
be a continuous function (often called a \emph{potential}).
Then
$$
P_T(\phi)=\sup_{\mu\in M_T(X)}\left(h_\mu(T)+\int \phi\,d\mu\right),
$$
where $M_T(X)$ denotes the collection of $T$-invariant Borel probability measures. 
\end{thm}
In particular, the special case $\phi=0$ gives the Variational Principle for topological entropy mentioned above.

Soon afterwards, in a short and elegant proof, Misiurewicz \cite{MisiurewiczVP} gave 
a new proof of Walters's theorem that
has become the standard proof of the theorem.
In the case that the supremum in the Variational Principle is attained by a measure $\mu$,
the measure $\mu$ is called an \emph{equilibrium measure}, borrowing the term from statistical physics. 

In the case of subshifts on finite alphabets, the map $\mu\mapsto h_\mu(T)$ is upper semi-continuous with respect to the 
weak$^*$ topology, and by definition, $\mu\mapsto \int\phi\,d\mu(T)$ is continuous. Since the collection of $T$-invariant
measures is weak$^*$-compact, the supremum in the variational principle
is attained by some invariant measure. More generally, this holds if the map $T$ has the 
\emph{expansive} property \cite{BowenEntropyExpansive}. 
In the case where there is an equilibrium measure, questions of uniqueness play an important role. 
A well known case in which there is a unique equilibrium measure is the case where $(X,T)$ is an
irreducible shift of finite type and $\phi$ is a H\"older continuous potential (\cite{RuelleThermo}).

\subsection{A simple heuristic motivating the Variational Principle}
In this section, we work with full shift spaces, and give a simple heuristic derivation of the Variational Principle.
By the binomial theorem, if $a$ and $b$ are positive numbers, then
$$
(a+b)^n=\sum_{j=0}^n \binom nj a^jb^{n-j};
$$
In fact, the significant contributions are tightly concentrated around the values where the 
largest value of the summand occurs,
$j\approx \frac a{a+b}n$ (one can check the maximum occurs here by taking the ratio of consecutive terms and setting it equal to 1).

Since the ratio between a sum over $n+1$ terms and the largest term of the sum is between 1 and $n+1$; and $\log (n+1)/n\to 0$, we see
\begin{equation}\label{eq:obvious}
\begin{split}
&\lim_{n\to\infty}\frac 1n\log\sum_{j=0}^n\binom nja^jb^{n-j}
=\lim_{n\to\infty}\frac 1n\max_{\alpha\in[0,1]}\log\left(\binom n{\alpha n} a^{\alpha n}b^{(1-\alpha)n}\right)\\
&=\lim_{n\to\infty}\frac 1n \log\left(\left(\left(\frac a{a+b}\right)^{-\frac {a}{a+b}}
\left(\frac {b}{a+b}\right)^{-\frac {b}{a+b}}\right)^na^{\frac a{a+b}n}b^{\frac{b}{a+b}n}\right)\\
&=\left(-\frac{a}{a+b}\log\frac{a}{a+b}-\frac b{a+b}\log\frac{b}{a+b}\right)+\left(\frac a{a+b}\log a+\frac b{a+b}\log b\right)\\
&=\log(a+b),
\end{split}
\end{equation}
as expected, where we made use of the rough approximation 
$\binom n{\alpha n}\approx(\alpha^{-\alpha}(1-\alpha)^{-(1-\alpha)})^n$, which can be 
verified from Stirling's formula. 

We now reinterpret this in terms of dynamical systems, let $T$ be the shift map on 
$\{0,1\}^\Z$ and $\phi$ be the function given by
$$
\phi(x)=\begin{cases}
\log a&\text{if $x_0=0$; and}\\
\log b&\text{if $x_0=1$.}
\end{cases}
$$
Using \eqref{eq:Pcalc}, we see $P_n(T,\phi,1)=\sum_{w\in\{a,b\}^n}\prod_{i=0}^{n-1} w_i=(a+b)^n$, so 
that $P_T(\phi)=\log(a+b)$, the last line of \eqref{eq:obvious}. We can also interpret the second last line
of \eqref{eq:obvious} as $h_T(\mu)+\int\phi\,d\mu$, where $\mu$ is the product measure giving
mass $\frac a{a+b}$ to the $0$ symbol and $\frac b{a+b}$ to the 1 symbol.

That is, the terms giving the largest contribution to $P_n(T,\phi,1)$ for large $n$ correspond 
to $\mu$ (which is the equilibrium state);
and the pressure is just the growth rate of these corresponding terms. 

More generally, staying with a shift map, if $\phi$ is a general continuous function and $\mu$ is an ergodic measure,
by the Shannon-MacMillan-Breiman theorem, the majority of the measure is supported on approximately $e^{h(\mu)n}$ 
$n$-cylinders (each with measure approximately $e^{-hn})$. By the Birkhoff ergodic theorem, the majority of the measure
is supported on cylinder sets with $S_n\phi(x)\approx n\int \phi\,d\mu$. The contribution to 
$P_n(T,\phi,1)$ from the $\mu$-typical
cylinders is therefore approximately $e^{nh(\mu)}\times e^{n\int\phi\,d\mu}=e^{n(h(\mu)+\int\phi\,d\mu)}$. 
It is not hard to make this argument rigorous, proving $P(\phi)\ge h(\mu)+\int\phi\,d\mu$. 
Since this holds for all ergodic invariant 
measures $\mu$, we obtain $P(\phi)\ge \sup_\mu (h(\mu)+\int\phi\,d\mu)$. From the fact that $\frac 1n
\log\sum_{i=1}^k c_ie^{a_i n}\to\max a_i$, one may hope for an equality here; 
this is what is established by the Variational Principle.
In terms of the heuristic, it suggests that the growth rate in the pressure sum over all 
cylinder sets is governed by the growth rate
in the pressure sum coming only from the cylinder sets that have substantial measure with respect to an equilibrium state.

\section{The article of Ledrappier and Walters}

The article of Ledrappier and Walters deals with the relative scenario. To introduce this, we briefly 
summarize the previous section. 

In the theory described so far, one starts with a continuous dynamical system $T\colon X\to X$. 
Given an invariant measure measure $\mu$, we defined $h_\mu(T)$ and the 
related concept of topological entropy, $h_\text{top}(T)$. 
The Variational Principle then shows that $h_\text{top}(T)$ is the supremum of the $h_\mu(T)$ over the 
collection of invariant measures. 
Analogously, given a potential $\phi$, one could define the ``pressure of a measure" by 
$h_\mu(T)+\int\phi\,d\mu$ and a definition is given for the (topological) pressure. Again, the Variational Principle
establishes that the topological pressure is the supremum of the pressures of invariant measures. 
Measures of maximal entropy and equilibrium states, that is those measures achieving the supremum, 
have arisen in a wide array of situations as exactly the measures of interest. 
Examples of this include information theory (Shannon coding), smooth dynamics (physical 
measures), mathematical physics (Gibbs measures).

In the relative situation, one considers a factor map $\pi$ between a continuous dynamical system 
$T\colon X\to X$ and a second system $S\colon Y\to Y$. That is, $\pi\colon X\to Y$ is continuous and onto
and satisfies $\pi\circ T=S\circ\pi$. Given an $S$-invariant measure, $\nu$ on $Y$, one may look at the
family of invariant measures on $X$ whose push-forward is $\nu$. This is non-empty and is typically a large set.
Analogous to the measure-theoretic entropy and pressure of a measure mentioned above, one
can study the relative entropy of $\mu$ over $\nu$ and the relative pressure of $\mu$ over $\nu$. 
As equilibrium states and measures of maximal entropy are important in the non-relative case,
the relative thermodynamic formlism provides a way to identify measures of interest in the
relative setting in which an invariant measure is specified for the factor transformation. 

Ledrappier and Walters prove the following extension of the Variational Principle for
topological entropy:

\begin{thm}[Ledrappier and Walters \cite{LedrappierWalters}, relative topological entropy version]
\label{th:LW-topological}
Let $T\colon X\to X$ and $S\colon Y\to Y$ be continuous maps of compact metric spaces. 
Suppose further that $\pi$ is a continuous factor map from $(X,T)$ to $(Y,S)$. 
If $\nu$ is an $S$-invariant measure, then we have
$$
\sup_{\pi_*\mu=\nu}h_\mu(T)=h_\nu(S)+\int_Y h_\text{top}(T,\pi^{-1}(y))\,d\nu(y).
$$
\end{thm}

A rough version of this could be ``The maximal entropy of a lift of $\nu$ is the entropy 
of $\nu$ plus the average fibre entropy".

A measure $\mu$ for which the supremum is achieved is called a 
\emph{relative measure of maximal entropy over $\nu$}.
They also established the following generalization to pressure.

\begin{thm}[\cite{LedrappierWalters}, relative pressure version]\label{thm:LW-pressure}
Let $T\colon X\to X$, $S\colon Y\to Y$, $\pi:X\to Y$ be as before
and suppose that $\phi\colon X\to\R$ is continuous and $\nu$ is an $S$-invariant measure. Then
$$
\sup_{\pi_*\mu=\nu} \left(h_\mu(T)+\int \phi\,d\mu\right)=h_\nu(S)+ \int_Y P_T(\phi,\pi^{-1}(y))\,d\nu(y).
$$
\end{thm}

A measure for which equality is obtained is called a \emph{relative equilibrium state for $\phi$ over $\nu$}.

These theorems generalize the standard Variational Principle for topological entropy 
and pressure respectively (taking $Y$ to be the one point space). 
Theorem \ref{th:LW-topological} also extends the following intuitive-sounding theorem of Bowen:

\begin{thm}[Bowen \cite{Bowen:EntropyforGroupEndomorphisms}]
Let $(X,d)$ and $(Y,d')$ be compact metric spaces and $T\colon X\to X$ and $S\colon Y\to Y$ 
be continuous. Let $\pi\colon X\to Y$ be a continuous factor map.
Then
$$
h_\text{top}(T)\le h_\text{top}(S)+\sup_{y\in Y}h_\text{top}(T,\pi^{-1}y).
$$
\end{thm}

\section{Subsequent developments}
The following subsections point to some areas where the relative variational principle has had a strong impact. 

\subsection{Extension entropies, fiber entropies and entropy structures}
 
Recall in the definitions of topological entropy and pressure, one studies the growth rate for $(n,\epsilon)$-spanning
sets for a fixed $\epsilon>0$, taking the limit as $n\to\infty$; and one then takes the limit as $\epsilon$ is shrunk to 0. 
This is sometimes referred to as ``dynamics at different scales": the amount of complexity seen depends on the separation
parameter. On the other hand, as described above, in the symbolic setting where $T$ is a shift map (and more generally
if $T$ satisfies the expansiveness property), the growth rate is 
the same for all $\epsilon$ sufficiently small. 

Early indications of the subtlety of this appear in work of Bowen 
\cite{BowenEntropyExpansive}
and Misiurewicz \cite{MisiurewicznotHExpansive,MisiurewiczConditional}.,
where one studies quantities of the form $h_\text{top}(T,\{y\colon d(T^nx,T^ny)<\epsilon\})$. 
That is, one is asking about the entropy of the set of points that $\epsilon$-shadow $x$ forever. 
If this is zero for all sufficiently small $\epsilon$, the transformation is said to be \emph{$h$-expansive} (a weakening
of the expansiveness condition),
and if it converges to 0, the transformation is \emph{asymptotically $h$-expansive}.
Asymptotic $h$-expansiveness is a sufficient condition for upper semi-continuity 
of the map $\mu\mapsto h_\mu(T)$ \cite{MisiurewiczConditional},
so guarantees the existence of measures of maximal entropy and equilibrium states.
Subsequent appearances of the phenomenon in a smooth setting are in Buzzi \cite{Buzzi-Intrinsic} and 
Newhouse \cite{Newhouse}, where it is shown that $C^\infty$ maps are asymptotically $h$-expansive. 

The situation for maps maps that fail to be $h$-expansive is quite subtle. 
A general abstract framework capturing this phenomenon is the entropy structures 
introducted by Downarowicz \cite{DownarowiczEntropyStructure}, and further studied in
 \cite{BurguetMcGoff} and \cite{Burguet}.
Entropy structures are built on the study of the symbolic extension entropy of a dynamical system. 
Given a topological dynamical system $T\colon X\to X$, one is interested in the 
subshifts $S\colon Y\to Y$ that factor continuously onto $T\colon X\to X$:
$$
h_\text{sex}(T)=\inf \{h_\text{top}(S)\colon (Y,S)\text{ is a symbolic extension of $(X,T)$}\},
$$
See \cite{Downarowicz2001, BoyleFiebig2,BoyleDownarowicz} for more information. 

Downarowicz and Newhouse \cite{DownarowiczNewhouse} 
and Burguet \cite{BurguetSymbExt} studied symbolic extensions
of smooth systems. 
\subsection{Partially hyperbolic dynamics}

The relative variational principle has been broadly applied in
the thermodynamic formalism of partially hyperbolic dynamical systems
\cite{Ures, BuzziFisher,BuzziFisherSambarinoVasquez,TahzibiYang,RHRHTU}.
Given a partially hyperbolic map, one may wish to quotient out by the center manifold and relate the dynamics
of the original system to the quotient system. The formulation of the relative variational principle means that 
it is well suited for questions of this type.

\subsection{Geometric measure theory and fractal dimensions}

The Relative Variational Principle turns out to have a highly satisfying connection to some problems in 
fractal geometry. Here, some classical objects of study are the Hausdorff dimension and the Minkowski (box) dimension. 
In many classical cases, such as the middle-third Cantor set, the two dimensions agree, and one can show
from the definition that the Hausdorff dimension is bounded above by the box dimension. 
We focus on the setting of self-affine fractals (that is fractal sets 
where there are a number of essentially disjoint pieces that are affine images of the entire fractal).
Of course, the middle-third Cantor set fits into this framework since it it is the disjoint union of two pieces, each
of which is a $\frac 13$-scaled copy of the whole.
Bedford \cite{Bedford} and McMullen \cite{McMullen} considered a generalization of this where the fractal is obtained
by iteratively contracting the unit square by the matrix $\operatorname{diag}(m^{-1},n^{-1})$ with $n\ge m$
and translating by a family of ``digits",
$D\subset\{(\frac im,\frac jn):0\le i<m;\ 0\le j<n\}$. They calculated the (generally distinct) box and Hausdorff
dimensions for these sets. 
Subsequently, Lalley and Gatzouras \cite{LalleyGatzouras} considered a generalized situation, again with 
stronger contraction in one coordinate direction than the other. Gatzouras and Peres \cite{GatzourasPeresSurvey}
were concerned with repellers for expanding dynamical systems (the Cantor sets mentioned above 
fit into this framework). They studied conditions for existence and uniqueness of invariant measures
supported on the repeller for which the Hausdorff dimension of the measure agrees with the 
Hausdorff dimension of the repeller. They pointed out a relationship with 
the work of Ledrappier and Young \cite{LedrappierYoung} and gave a
reformulation of their central question in the form:
\begin{problem*}
Let $(X,T)$ be a subshift and let $(Y,S)$ be a factor by a map $\pi$. Let $\phi\colon X\to\R$ and
$\psi:Y\to\R$ satisfy $\phi(x)\ge \psi(\pi(x))>0$. Is the supremum
$$
\sup_{\mu\in M_T(X)}\left(\frac{h_\mu}{\int\phi\,d\mu}+\left(\frac 1{\int\psi\,d\pi_*\mu}-\frac 1{\int\phi\,d\mu}\right)
h_{\pi_*\mu}\right)
$$
attained for some ergodic measure $\mu$?
\end{problem*}
In the application of statement to the calculation of Hausdorff dimensions,
the $\phi(x)$ and $\psi(\pi(x))$ should be thought of as partial derivatives in
the coordinate directions, and $\pi$ is the projection onto the first coordinate. 

Shin \cite{ShinETDS} studied the special case where $\phi$ is constant, giving sufficient conditions
in terms of compensation functions (see the next subsection) for
a positive answer. The question was further studied in \cite{Feng}, giving a positive answer in terms
of specification conditions on $X$. Feng and Huang \cite{FengHuang} extended the Relative Variational
Principle to what they called \emph{weighted pressure}, a generalization of the above,
and used the extension to compute Hausdorff dimension of some affine repellers. 

\subsection{Compensation functions and relative equilibrium states}
Walters \cite{WaltersCompensation} studied a notion of \emph{compensation function} that is
related to Theorem \ref{thm:LW-pressure} and a concept that was first studied in 
\cite{BoyleTuncel}. 
In the context of Theorem \ref{thm:LW-pressure}, a compensation is a function $\phi$ satisfying
$P_T(f\circ\pi+\phi)=P_S(f)$ for all $f\in C(Y)$. That is, $\phi$ ``compensates" for the (possible) additional 
entropy in the preimage fibres. These compensation functions take a particularly nice form 
in the context of \cite{BoyleTuncel}. Measurable compensation functions always
exist as a corollary of Theorem \ref{thm:LW-pressure}. As is common in thermodynamic formalism,
it is of interest to look for compensation functions with more regularity or other features making them easier
to work with. 
Antonioli \cite{AntonioliCompensation} studied the regularity
of the compensation functions in \cite{WaltersCompensation} and Shin \cite{Shin} studied saturated compensation functions.

In the non-relative setting, the multiplicity of equilibrium states
has been studied by many authors. The corresponding question in the relative setting
was addressed in \cite{PetersenShin,PQS,AllahbakhshiQuas,AllahbakhshiAntonioliYoo,YooJMD,YooTAMS}.

\bibliographystyle{abbrv}
\bibliography{LedrappierWaltersbib}
\end{document}